\documentclass[11pt,a4paper]{amsart}

\usepackage[utf8]{inputenc}
\usepackage[T1]{fontenc}
\usepackage{amsmath,amssymb,amsthm,mathtools}
\usepackage[margin=2.5cm]{geometry}
\usepackage{enumitem}
\usepackage{booktabs,tabularx}
\usepackage{float} 
\usepackage{microtype}
\usepackage{xcolor} 
\usepackage{tikz}
\usepackage[colorlinks=true,linkcolor=blue!50!black,citecolor=green!40!black,urlcolor=blue!60!black]{hyperref}

\newtheorem{theorem}{Theorem}[section]
\newtheorem{proposition}[theorem]{Proposition}
\newtheorem{lemma}[theorem]{Lemma}
\newtheorem{corollary}[theorem]{Corollary}
\theoremstyle{definition}

\newtheorem{remark}[theorem]{Remark}
\newtheorem{example}[theorem]{Example}
\numberwithin{equation}{section}

\newcommand{\R}{\mathbb{R}}
\newcommand{\Z}{\mathbb{Z}}

\newcommand{\T}{\mathbb{T}}
\newcommand{\dd}{\,\mathrm{d}}
\newcommand{\E}{\mathcal{E}}
\newcommand{\X}{\mathcal{X}}
\newcommand{\Rset}{\mathcal{R}(T)}
\newcommand{\Out}{\mathcal{O}}
\newcommand{\supp}{\operatorname{supp}}
\newcommand{\norm}[1]{\left\|#1\right\|}
\newcommand{\dist}[1]{\| #1 \|_{\Z}}

\title[Reachable sets of point sources on the line]{Controllability of the 1D wave equation by point sources: outgoing states, circle rotations and small divisors}

\author{Kaïs Ammari}
\address{LR Analysis and Control of PDEs, LR22ES03, Department of Mathematics, Faculty of Sciences of Monastir, University of Monastir, Tunisia}
\email{kais.ammari@fsm.rnu.tn}

\author{Mourad Sini}
\address{RICAM, Austrian Academy of Sciences, Altenbergerstrasse 69, 4040, Linz, Austria}
\email{mourad.sini@oeaw.ac.at}

\date{November 15, 2026}


\newcommand{\AT}{\mathcal A_T}
\newcommand{\YT}{\mathcal Y_T}
\newcommand{\GT}{\mathcal G_T}
\newcommand{\JT}{\mathcal J_T}
\newcommand{\CT}{\mathfrak C_T}
\newcommand{\DK}{\mathfrak D}
\newcommand{\Kres}{K_{\mathrm{res}}}
\newcommand{\ran}{\operatorname{Ran}}
\newcommand{\ind}{\mathbf 1}
\newcommand{\br}[1]{\langle #1\rangle}

\newcommand{\Q}{\mathbb Q}
\newcommand{\id}{\operatorname{Id}}
\hypersetup{
  pdftitle={Controllability of the 1D wave equation by point sources: outgoing states, circle rotations and small divisors},
  pdfauthor={Kais Ammari and Mourad Sini},
  pdfsubject={Reachability, circle rotations, small divisors, and quantitative control cost},
  pdfkeywords={wave equation, point sources, controllability, small divisors, Diophantine approximation}
}
\allowdisplaybreaks[2]
\begin{document}

\begin{abstract}
We study the terminal map $S_Tq=(w(T),w_t(T))$ for
\[
k^{-1}w_{tt}-\rho^{-1}w_{xx}
 =\sum_{j=1}^M q_j(t)\delta_{z_j}
 \quad\text{in }(0,\infty)\times\R,
 \qquad w(0)=w_t(0)=0,
\]
where $k,\rho>0$, $z_1<\cdots<z_M$, $M\ge2$, and $q\in L^2(0,T)^M$ is extended by zero. Set $c=\sqrt{k/\rho}$, $d=z_M-z_1$, $d_j=z_j-z_1$, and $T_*=2d/c$. We characterize $\Rset=\ran S_T$ in the compatible outgoing-state space $\AT\subset\X=\dot H^1(\R)\times L^2(\R)$, defined by the finite-propagation supports of the Riemann invariants and $\int_\R u_0'(x)\dd x=0$. Our main results are as follows.
\begin{enumerate}[label=\textup{(\roman*)},leftmargin=*,itemsep=3pt,topsep=4pt]
\item Let $\Pi_U$ be the $2d$-periodisation of $R(z_1+\cdot)-L(z_1-\cdot)$, where $R,L=\tfrac12(u_1\mp cu_0')$. A state $U\in\AT$ is reachable if and only if
\[
\Pi_U(y)=\sum_{j=2}^{M-1}\bigl(f_j(y-2d_j)-f_j(y)\bigr),
\qquad f_j\in L^2(\T_d),\quad
\supp f_j\subset[-d_j,cT-d_j]\pmod{2d},
\]
where $\T_d=\R/(2d\Z)$. The controls are recovered constructively from the $f_j$ and a bounded two-source reconstruction.
\item For $M=2$, $S_T$ is an isomorphism onto its closed range for every $T>0$. For $M=3$, the range is closed when $T<T_*$. For arbitrary $M\ge3$ and $T\ge T_*$, it is closed if and only if all $\theta_j:=d_j/d$ are rational, and it is dense in $\AT$ if and only if at least one $\theta_j$ is irrational.
\item For $M\ge3$, $T\ge T_*$, and $U\in\Rset$, the minimum control norm satisfies
\[
\CT(U)\asymp\norm U_\X+
\left(2d\sum_{\sigma_n>0}\frac{|\widehat\Pi_{U,n}|^2}{\sigma_n^2}\right)^{1/2},
\qquad
\sigma_n^2=4\sum_{j=2}^{M-1}\sin^2(\pi n\theta_j).
\]
\item Put $p=M-2$. In the same time regime, the simultaneous condition
$\max_j\dist{n\theta_j}\ge b|n|^{-1/p}$ for all $n\in\Z\setminus\{0\}$ and some $b>0$ gives a sharp worst synthesis cost of order $N^{1/p}$ on energy-normalised states with defects of Fourier degree at most $N$. On the canonical states associated with the unit $H^s$ ball of mean-zero defects, the optimal cost of approximation to energy error $\varepsilon$ is of order $\varepsilon^{-(1/p-s)/s}$ for $0<s<1/p$ and the threshold $s=1/p$ for uniformly bounded exact synthesis is sharp.
\end{enumerate}
\end{abstract}

\maketitle

\section{Introduction}

\subsection{Mathematical setting}
Let $\rho,k>0$, $c=\sqrt{k/\rho}$, and let $z_1<\dots<z_M$ be points of the real line. For controls $q=(q_1,\dots,q_M)\in L^2(0,T)^M$, extended by zero outside $(0,T)$, we consider
\begin{equation}\label{eq:main}
\begin{cases}
k^{-1}w_{tt}-\rho^{-1}w_{xx}=\displaystyle\sum_{j=1}^Mq_j(t)\,\delta_{z_j},
       & (t,x)\in(0,\infty)\times\R,\\
w(0,\cdot)=w_t(0,\cdot)=0.
\end{cases}
\end{equation}
We set $S_Tq=(w(T),w_t(T))$ and denote its range by
\[
\Rset=\{S_Tq:q\in L^2(0,T)^M\}.
\]
We seek an intrinsic criterion for membership in $\Rset$, a classification of its closedness and of the uniqueness of the controls, and quantitative estimates for the cost of reaching or approximating a prescribed state.

The energy space is $\X=\dot H^1(\R)\times L^2(\R)$, where $\dot H^1(\R)$ consists of locally $H^1$ functions with derivative in $L^2$, modulo additive constants. We use
\begin{equation}\label{eq:energy}
\E(U)=\frac12\int_\R\big(k^{-1}|u_1|^2+\rho^{-1}|u_0'|^2\big)\dd x,
\qquad \norm U_\X^2=2\E(U),\qquad U=(u_0,u_1).
\end{equation}
All states and controls may be taken real; complex Fourier coefficients are used only for convenience. For $M\ge2$, write
\[
 a_{jl}=\frac{|z_j-z_l|}{c},\qquad d=z_M-z_1,\qquad d_j=z_j-z_1,
 \qquad \theta_j=\frac{d_j}{d}=\frac{a_{1j}}{a_{1M}}\quad(1<j<M).
\]
The case $M=1$, which does not require this notation, is treated in Section~\ref{sec:outgoing}.

For $U\in\X$, introduce the Riemann invariants
\begin{equation}\label{eq:invariants}
R=\tfrac12(u_1-cu_0'),\qquad L=\tfrac12(u_1+cu_0'), \text{ then }
\quad \norm U_\X^2=\frac2k\big(\norm R_2^2+\norm L_2^2\big).
\end{equation}
Every reachable state satisfies
\begin{equation}\tag{$\mathrm S_T$}\label{eq:S}
\supp R\subset[z_1,z_M+cT],\qquad
\supp L\subset[z_1-cT,z_M].
\end{equation}
Let $\YT$ be the closed subspace of $\X$ defined by \eqref{eq:S}. The correct ambient space for the density statements below is
\begin{equation}\label{eq:AT}
\AT=\left\{U\in\YT:\int_\R u_0'(x)\dd x=0\right\}
     =\left\{U\in\YT:\int_\R R\dd x=\int_\R L\dd x\right\}.
\end{equation}
The integrals are well defined because the supports are bounded. The additional condition is necessary as it says that the displacement has a representative which vanishes at both spatial infinities. In particular, $\AT$ is closed in $\X$, and its elements have a unique compactly supported displacement representative. The support condition alone does not imply this compatibility in the homogeneous energy space.

\subsection{Background}
Pointwise wave control is a classical instance in which actuator geometry and arithmetic interact. For a finite string with Dirichlet endpoints, an actuator at $\xi$ couples to the $n$th eigenfunction through $\sin(n\pi\xi/\ell)$. Vanishing coefficients obstruct approximate controllability, whereas coefficients approaching zero affect exact reachability and the regularity of reachable states. The general theory of pointwise control and the moment method is developed in \cite{Lions88,Lions92,AvdoninIvanov,TucsnakWeiss} and related questions for one-dimensional strings are treated in \cite{FabrePuel,JTZ97}.

\noindent Two comparisons are particularly relevant here. Avdonin and Seidman \cite{AvdoninSeidman} study pointwise control as a limit of spatially distributed actuation and explain how the choice of the data spaces affects sensitivity to actuator location. Avdonin and Tucsnak \cite{AvdoninTucsnak01} determine sharp-time simultaneous controllability properties for two bounded strings driven from a common endpoint, including arithmetic conditions on the reachable subspace. These works show that the relation between arithmetic and reachability is not specific to the present model. Our concern is instead the complete reachable set for independent fixed point sources on the \emph{whole line}, where the reduction comes from cancellation of outgoing waves rather than boundary reflections or a bounded-domain eigenfunction expansion.

\noindent  Arithmetic also enters pointwise stabilisation, but the conclusions depend on the boundary conditions and the functional setting. For a string with mixed Dirichlet--Neumann boundary conditions, Ammari, Henrot and Tucsnak \cite{AmmariHenrotTucsnak} characterise actuator positions yielding exponential stability and obtain polynomial estimates in nonexponential regimes. Further results on singular damping and unbounded feedbacks are given in \cite{JTZ98,AmmariTucsnak01,AmmariNicaise}. Point masses, networks and moving actuators lead to related, but distinct, control problems; see \cite{HansenZuazua,DagerZuazua,Khapalov95,Castro13}.

\noindent  The whole-line problem has no reflecting boundary. Nevertheless, three sources already produce an arithmetic obstruction as  eliminating the outer controls leaves a translation equation on a circle, and the inner source determines its rotation number. This gives both a finite-time classification and a quantitative interpretation of small divisors in terms of the physical $L^2$-norm of the controls.

\subsection{Main results}
Each source produces a closed subspace $\Out_j(T)$ of even, compactly supported outgoing states, and
\[
\Rset=\Out_1(T)+\dots+\Out_M(T).
\]
The question is whether a prescribed state admits such a decomposition and how stably that decomposition can be obtained.

For $U\in\YT$, define its defect and its periodisation on $\T_d=\R/(2d\Z)$ by
\begin{equation}\label{eq:defect}
\Delta_U(y)=R(z_1+y)-L(z_1-y),\qquad
\Pi_U(y)=\sum_{n\in\Z}\Delta_U(y+2nd).
\end{equation}
We often omit the subscript $U$. The defect compares the right-going part of a state with the mirror image of its left-going part about $z_1$. The map $\GT:U\mapsto\Pi_U$ is bounded on $\YT$, and $\Pi_U$ has mean zero precisely when $U\in\AT$. For an inner source, let
\[
A_j(T)=[-d_j,cT-d_j]\pmod{2d};
\]
this is the whole circle when $cT\ge2d$.

\begin{theorem}[Characterisation of $\Rset$]\label{thm:main}
Let $M\ge2$. A state $U\in\AT$ belongs to $\Rset$ if and only if there exist $f_j\in L^2(\T_d)$, $1<j<M$, vanishing outside $A_j(T)$, such that
\begin{equation}\label{eq:cohom}
\Pi(y)=\sum_{j=2}^{M-1}\big(f_j(y-2d_j)-f_j(y)\big)
\quad\text{for a.e. }y\in\T_d.
\end{equation}
No state outside $\AT$ is reachable. Given such $f_j$, put $\ell_T=\min(cT,2d)$. One choice of the inner controls is
\begin{equation}\label{eq:inner-control}
q_j(s)=\frac{2}{\rho c}\,
\ind_{(0,\ell_T)}\big(c(T-s)\big)
 f_j\big(c(T-s)-d_j\big),\qquad 0<s<T,
\end{equation}
where $f_j$ is evaluated periodically. Subtracting the resulting inner-source states, the outer controls are recovered from \eqref{eq:recursion}, applied to the pair $(z_1,z_M)$.
\end{theorem}

The support restriction in Theorem~\ref{thm:main} is essential for short
control times. At $T\ge 2d/c$ it disappears, and the remaining obstruction
can be described by Fourier multipliers. The other results are summarized
below.
\begin{enumerate}[label=\textup{(\roman*)},leftmargin=*,itemsep=5pt,topsep=5pt]
\item \emph{Two sources and the defect quotient.}
For $M=2$, the reachable set is $\ker\GT$, is closed in $\AT$, and every
reachable state admits a unique control satisfying
$\|q\|_{L^2(0,T)^2}\le C_T\|U\|_{\X}$ for all $T>0$
(Theorem~\ref{thm:two}). More generally, letting
$\mathcal R_{\mathrm{out}}(T)=\Out_1(T)+\Out_M(T)$, we prove
$\mathcal R_{\mathrm{out}}(T)=\ker(\GT|_{\AT})$ and, for $T\ge 2d/c$,
\[
\AT/\mathcal R_{\mathrm{out}}(T)\simeq L^2_0(\T_d),
\]
with a bounded right inverse $\JT$ of $\GT$
(Lemma~\ref{lem:quotient}).

\item \emph{Three sources at arbitrary control times.}
For $M=3$, the support constraint in~\eqref{eq:cohom} gives a
finite-orbit criterion below $T_*=2a_{13}$.
Table~\ref{tab:three} states the complete range and uniqueness classification,
including the critical times (Theorem~\ref{thm:three}); here
$\theta=a_{12}/a_{13}$ and, for rational $\theta=m/v$ in lowest terms,
$\tau=a_{13}/v$.

\begin{table}[H]
\centering\small
\renewcommand{\arraystretch}{1.20}
\begin{tabularx}{\textwidth}{@{}p{2.9cm}p{3.8cm}X@{}}
\toprule
Rotation number & Time range & Reachable set and controls\\
\midrule
$\theta\notin\Q$ & $T<2a_{13}$ & Closed; unique controls.\\
 & $T\ge2a_{13}$ & Dense in $\AT$, not closed; non-unique controls.\\
\midrule
$\theta=m/v$ & $T\le2(a_{13}-\tau)$ & Closed; unique controls.\\
(in lowest terms) & $2(a_{13}-\tau)<T<2a_{13}$ & Closed; non-unique controls.\\
 & $T\ge2a_{13}$ & Closed, proper in $\AT$; non-unique controls.\\
\bottomrule
\end{tabularx}
\caption{Three-source classification; $\tau=a_{13}/v$ in the rational case. Closedness is in the energy topology.}
\label{tab:three}
\end{table}

\item \emph{Long-time reachability and the physical control cost.}
For $M\ge3$ and $T\ge 2d/c$, set
$\theta_j=d_j/d$ and
$\sigma_n^2=4\sum_{j=2}^{M-1}\sin^2(\pi n\theta_j)$.
A state $U\in\AT$ is reachable precisely when
$\widehat\Pi_{U,n}=0$ at every common resonance $\sigma_n=0$ and
$\sum_{\sigma_n>0}|\widehat\Pi_{U,n}|^2/\sigma_n^2<\infty$.
The range is closed exactly when all $\theta_j$ are rational, and is dense
in $\AT$ exactly when at least one $\theta_j$ is irrational
(Theorem~\ref{thm:long}). For reachable $U$, define
\begin{equation}\label{eq:control-cost}
\CT(U)=\min\{\norm q_{L^2(0,T)^M}:S_Tq=U\}.
\end{equation}
Theorem~\ref{thm:physical-cost} establishes the two-sided estimate
\begin{equation}\label{eq:intro-cost}
\CT(U)\asymp \norm U_\X+
\left(2d\sum_{\sigma_n>0}
\frac{|\widehat\Pi_{U,n}|^2}{\sigma_n^2}\right)^{1/2},
\end{equation}
with constants depending on the fixed medium, time and configuration.

\item \emph{Sharp arithmetic cost--accuracy estimates.}
Let $p=M-2$. Under the simultaneous badly approximable condition
$\max_j\dist{n\theta_j}\ge b|n|^{-1/p}$ for $n\ne0$, the
maximal energy-normalized cost for a defect of Fourier degree at most $N$
is of order $N^{1/p}$ (Proposition~\ref{prop:multi-arithmetic}).
For $0<s<1/p$, the optimal approximation cost on the canonical unit
$H^s(\T_d)$ ball satisfies
\[
\sup_{\substack{g\in H^s(\T_d),\ \widehat g_0=0\\
                  \norm g_{H^s}\le1}}
\ \inf_{\norm{S_Tq-\JT g}_\X\le\varepsilon}
\norm q_{L^2(0,T)^M}
\asymp \varepsilon^{-(1/p-s)/s}.
\]
Here the lower bound holds for \emph{all} controls, and even without the
badly approximable condition provided some $\theta_j$ is irrational.
The threshold $s=1/p$ for uniformly bounded exact synthesis is sharp and
for three sources, Liouville positions have superpolynomial approximation
costs along subsequences (Theorem~\ref{thm:minimax}).
\end{enumerate}

\medskip
The arithmetic dependence of controllability and stabilization on actuator locations is already known for bounded strings (see Section~1.2). Our setting differs in that the quotient by the closed outer-source range reduces reachability of the full terminal energy state to a supported cohomological equation on the circle. Moreover, \eqref{eq:intro-cost} quantifies small-divisor losses in terms of the \emph{minimum norm of the original point-source controls}, including the outer correction, and the sharp approximation lower bounds apply to all admissible controls. This one-dimensional structure contrasts with the three-dimensional
problem studied by Ammari and Sini in~\cite{AmmariSini3D}, where, for fixed distinct sources, $S_T^{(3)}:H_0^1(0,T)^M\to\X_3$ satisfies $\|S_T^{(3)}q\|_{\X_3}\asymp\|q'\|_{L^2(0,T)^M}$ and hence has an injective, closed range. In one dimension, the corresponding map on $L^2(0,T)^M$ may instead have a nontrivial kernel and a dense, nonclosed range. Thus the distinction involves both the reachable-state structure and the topology of the control-to-state map.

\subsection*{Outline}
Section~\ref{sec:outgoing} describes the outgoing subspaces and the case of short times. Section~\ref{sec:two} treats two sources, including stable reconstruction. Section~\ref{sec:main} proves Theorem~\ref{thm:main} and identifies the defect quotient of the state space. Section~\ref{sec:three} gives the three-source classification and the scalar arithmetic estimates. Section~\ref{sec:long} treats arbitrary numbers of sources at long times, proves the comparison with physical control cost, and establishes constructive and sharp approximation bounds.

\section{Outgoing states}\label{sec:outgoing}

\subsection{The control-to-state map}
The forward fundamental solution of $k^{-1}(\partial_t^2-c^2\partial_x^2)$ is $\frac{\rho c}{2}H(ct-|x|)$ where $H$ is the Heaviside function. With
\[
Q_j(t)=\int_0^t q_j(s)\dd s\quad(t\ge0),\qquad Q_j(t)=0\quad(t<0),
\]
the solution of \eqref{eq:main} is
\begin{equation}\label{eq:solution}
w(t,x)=\frac{\rho c}{2}\sum_{j=1}^M Q_j\left(t-\frac{|x-z_j|}{c}\right).
\end{equation}
The formula defines the distributional solution with zero initial data and its spatial derivative and time derivative belong to $L^2$ at each finite time. Differentiation, first for smooth controls and then by $L^2$ approximation, gives
\begin{equation}\label{eq:RL}
\begin{aligned}
R(x)&=\frac{\rho c}{2}\sum_{j=1}^M q_j\left(T-\frac{x-z_j}{c}\right),\\
L(x)&=\frac{\rho c}{2}\sum_{j=1}^M q_j\left(T+\frac{x-z_j}{c}\right).
\end{aligned}
\end{equation}
In the first line, terms with $x<z_j$ vanish because the argument exceeds $T$; in the second line, terms with $x>z_j$ vanish for the same reason. Formula \eqref{eq:solution} also shows that the displacement is supported in $[z_1-cT,z_M+cT]$. Consequently $S_Tq\in\AT$.

\begin{proposition}\label{prop:one}
The states produced by source $j$ alone form the closed subspace
\[
\begin{aligned}
\Out_j(T)=\Big\{(u_0,u_1):\;&u_1\in L^2(\R),\quad
\supp u_1\subset[z_j-cT,z_j+cT],\\
&u_1\text{ is even about }z_j,\qquad
u_0(x)=\frac1c\int_{|x-z_j|}^{cT}u_1(z_j+s)\dd s\Big\},
\end{aligned}
\]
where $u_0$ denotes the compactly supported representative. The control is recovered by
\[
q_j(s)=\frac2{\rho c}\,u_1\big(z_j+c(T-s)\big),\qquad 0<s<T,
\]
and
\begin{equation}\label{eq:one-energy}
\E(S_Tq_j)=\frac{\rho c}{2}\norm{q_j}_{L^2(0,T)}^2,
\qquad
\norm{S_Tq}_\X\le\sqrt{M\rho c}\,\norm q_{L^2(0,T)^M}.
\end{equation}
\end{proposition}

\begin{proof}
For a single source, \eqref{eq:RL} gives $L=0$ to the right of $z_j$, $R=0$ to its left, and $L(z_j-y)=R(z_j+y)$. Thus $u_1=-cu_0'$ on the right, $u_1=cu_0'$ on the left, and $u_1$ is even about $z_j$. Integrating from spatial infinity gives the formula for $u_0$. Conversely, the stated control produces each such state. A change of variables in \eqref{eq:RL} gives
\[
\norm R_2^2=\norm L_2^2=\frac{\rho^2c^3}{4}\norm{q_j}_2^2.
\]
Together with $k=\rho c^2$, this proves the energy identity, which also implies closedness of the one-source range. Finally,
\[
\norm{S_Tq}_\X\le\sum_j\norm{S_Tq_j}_\X
\le\sqrt{M\rho c}\left(\sum_j\norm{q_j}_2^2\right)^{1/2}.
\]
\end{proof}

\begin{remark}[Compact representatives and the energy topology]\label{rem:compatibility}
If $U\in\YT$, then $u_0'$ is supported in the fixed bounded interval $I_T=[z_1-cT,z_M+cT]$. A representative is constant on each component of $\R\setminus I_T$; the difference of these constants is $\int u_0'$. Hence \eqref{eq:AT} is equivalent to the existence of a compactly supported representative. The integral of $u_0'$ is a continuous functional on $\YT$, so $\AT$ is closed. On these compact representatives, Poincar\'e's inequality shows that the energy norm and the usual $H^1\times L^2$ norm are equivalent, with a constant depending on $I_T$. Without the zero-integral condition, this conclusion is false: $u_1=0$ and $u_0'=g$ with $g\in C_c^\infty(z_1,z_M)$ and $\int g\ne0$ satisfy \eqref{eq:S}, but their state is not in the closure of $\Rset$.
\end{remark}

\subsection{Short times}
\begin{corollary}\label{cor:short}
Suppose $M\ge2$ and $cT\le d_{\min}=\min_{j\ne l}|z_j-z_l|$. A compactly supported state belongs to $\Rset$ if and only if
\[
\begin{gathered}
\supp R\subset\bigcup_j[z_j,z_j+cT],\qquad
\supp L\subset\bigcup_j[z_j-cT,z_j],\\
L(z_j-y)=R(z_j+y)\quad\text{for a.e. }0<y<cT,\quad j=1,\dots,M.
\end{gathered}
\]
The sum $\Out_1(T)\oplus\dots\oplus\Out_M(T)$ is orthogonal in the energy inner product. In particular, the controls are unique and
\[
\norm{S_Tq}_\X^2=\rho c\sum_{j=1}^M\norm{q_j}_2^2.
\]
\end{corollary}

\begin{proof}
The intervals supporting the right-going invariants of distinct sources are disjoint up to endpoints; the same holds for the left-going invariants. Thus only source $j$ contributes to $R$ on $(z_j,z_j+cT)$ and to $L$ on $(z_j-cT,z_j)$. Apply Proposition~\ref{prop:one} on each pair of intervals. Orthogonality follows from \eqref{eq:invariants}; possible overlap between a right-going and a left-going wave does not create a cross term in that inner product.
\end{proof}

For larger times the right-going supports, and likewise the left-going supports, can overlap. The decomposition can then no longer be read off locally. This is the role of the defect.

\section{Two sources}\label{sec:two}

It is convenient to read the invariants as functions of time. Set
\begin{equation}\label{eq:time-invariants}
\begin{aligned}
r(\sigma)&=\frac2{\rho c}R\big(z_1+c(T-\sigma)\big),\\
l(\sigma)&=\frac2{\rho c}L\big(z_1+c(\sigma-T)\big),
\qquad t_j=\frac{d_j}{c}.
\end{aligned}
\end{equation}
Then \eqref{eq:RL} becomes
\begin{equation}\label{eq:rl}
r=\sum_jq_j(\cdot+t_j),\qquad l=\sum_jq_j(\cdot-t_j).
\end{equation}
Moreover,
\[
\Delta(y)=\frac{\rho c}{2}(r-l)(T-y/c),
\]
and \eqref{eq:S} becomes $\supp r\subset[-t_M,T]$ and $\supp l\subset[0,T+t_M]$.

\begin{theorem}[Two sources]\label{thm:two}
Let $M=2$ and $a=a_{12}$. A state $U\in\AT$ belongs to $\Rset$ if and only if $\Pi_U=0$, or equivalently
\[
\sum_{n\in\Z}R(z_1+y+2nd)=\sum_{n\in\Z}L(z_1-y+2nd)
\quad\text{for a.e. }y\in\R.
\]
The controls are unique and are given, as zero-extended functions on $\R$, by
\begin{equation}\label{eq:recursion}
q_2(s)=\sum_{n\ge0}(r-l)\big(s-(2n+1)a\big),
\qquad q_1=r-q_2(\cdot+a).
\end{equation}
The reachable set is closed for every $T>0$. In addition, the reconstruction is bounded: for a constant $B_T$ depending only on $T,a,\rho,c$,
\begin{equation}\label{eq:outer-bound}
\norm{(q_1,q_2)}_{L^2(0,T)^2}\le B_T\norm U_\X.
\end{equation}
\end{theorem}

\begin{proof}
\emph{Necessity.} By \eqref{eq:rl}, $D:=r-l=q_2(\cdot+a)-q_2(\cdot-a)$. Its $2a$-periodisation vanishes by cancellation of two translates. The change of variables above identifies this with $\Pi_U=0$.

\emph{Sufficiency.} Assume $\Pi_U=0$ and define \eqref{eq:recursion}. Since $D$ is supported in $[-a,T+a]$, the sum is locally finite. Telescoping gives
\[
q_2(s+a)-q_2(s-a)=D(s).
\]
Consequently $q_1+q_2(\cdot-a)=r-D=l$, so \eqref{eq:rl} holds. We check the control supports. If $s<0$, all arguments in \eqref{eq:recursion} lie below $-a$, and $q_2(s)=0$. If $s>T$, the omitted terms of the full $2a$-periodisation have arguments at least $s+a>T+a$, so the sum equals the full periodisation and vanishes. Therefore $\supp q_2\subset[0,T]$. For $s>T$, both terms defining $q_1$ vanish. On $(-a,0)$, the difference identity and $q_2(s-a)=0$ yield $q_1(s)=r(s)-D(s)=l(s)=0$. For $s<-a$ both terms vanish. Thus $q_1,q_2\in L^2(0,T)$, and their state has the prescribed invariants and hence is $U$ in $\X$.

\emph{Uniqueness.} The difference equation, together with $q_2=0$ on $(-\infty,0)$, forces the recursion \eqref{eq:recursion}; $q_1$ is then determined.

\emph{Bounded reconstruction and closedness.} On $0<s<T$, at most
\[
N_T^{\rm out}=1+\left\lfloor\frac{T}{2a}\right\rfloor
\]
terms can contribute to $q_2$. Hence
\[
\norm{q_2}_2\le N_T^{\rm out}\norm D_2,
\qquad \norm{q_1}_2\le\norm r_2+\norm{q_2}_2.
\]
The scaling \eqref{eq:time-invariants} gives
\[
\norm r_2+\norm l_2
\le\frac{2\sqrt{k}}{\rho c^{3/2}}\norm U_\X.
\]
Thus one may take
\[
B_T=\frac{2\sqrt{k}}{\rho c^{3/2}}(1+2N_T^{\rm out}).
\]
Finally, periodisation is bounded on functions supported in a fixed bounded interval: by Cauchy--Schwarz, its squared norm is at most the maximum overlap multiplicity times the squared norm of the function. Therefore $U\mapsto\Pi_U$ is continuous on $\AT$, and its kernel is closed. This is exactly the reachable set.
\end{proof}

The same theorem and estimate will be used for the outer pair $(z_1,z_M)$, with $a=a_{1M}$ and $q_2$ relabelled as $q_M$.

\section{Proof of the characterisation}\label{sec:main}

\subsection{Periodisation and the contribution of one source}
For later use set
\begin{equation}\label{eq:lift-constants}
\beta=\frac{2}{\rho c^{3/2}},\qquad
m_T=\left\lceil\frac{cT}{2d}\right\rceil\ge1.
\end{equation}

\begin{lemma}\label{lem:contrib}
Let
\[
\varphi_j(s)=\frac{\rho c}{2}q_j(T-s/c),\qquad
f_j(y)=\sum_{n\in\Z}\varphi_j(y+d_j+2nd).
\]
Then $\varphi_j$ is supported in $[0,cT]$, $f_j$ vanishes outside $A_j(T)$, and the periodised defect produced by source $j$ is
\begin{equation}\label{eq:one-defect}
\Pi^{(j)}(y)=f_j(y-2d_j)-f_j(y).
\end{equation}
In particular, $\Pi^{(1)}=\Pi^{(M)}=0$. Moreover,
\begin{equation}\label{eq:period-bound}
\norm{f_j}_{L^2(\T_d)}\le\sqrt{m_T}\norm{\varphi_j}_{L^2(\R)}
=\frac{\sqrt{m_T}}{\beta}\norm{q_j}_{L^2(0,T)}.
\end{equation}
Conversely, every $f_j\in L^2(\T_d)$ supported in $A_j(T)$ has the lift
\begin{equation}\label{eq:canonical-lift}
\varphi_j(s)=\ind_{(0,\ell_T)}(s)f_j(s-d_j),\qquad
\ell_T=\min(cT,2d),
\end{equation}
whose periodisation after translation by $d_j$ is $f_j$. This lift satisfies
\begin{equation}\label{eq:lift-isometry}
\norm{\varphi_j}_2=\norm{f_j}_{L^2(\T_d)},
\qquad \norm{q_j}_2=\beta\norm{f_j}_{L^2(\T_d)}.
\end{equation}
\end{lemma}

\begin{proof}
By \eqref{eq:RL}, the unperiodised defect of source $j$ is
$\varphi_j(y-d_j)-\varphi_j(y+d_j)$. Periodising and writing
$y-d_j=(y-2d_j)+d_j$ gives \eqref{eq:one-defect}. The two outer shifts are zero and one full period.

For a.e. $y$, at most $m_T$ translates of a function supported in an interval of length $cT$ contribute to its $2d$-periodisation. Cauchy--Schwarz followed by integration over one period proves \eqref{eq:period-bound}; the norm of $\varphi_j$ follows by changing variables from $s$ to $T-s/c$.

If $cT<2d$, the interval $(0,cT)-d_j$ maps injectively, up to endpoints, onto $A_j(T)$, and the support constraint proves the lifting assertion and the norm identity. If $cT\ge2d$, the interval $(0,2d)-d_j$ is one full period, so the same statements hold without a support constraint. Replacing $s$ by $c(T-t)$ gives precisely \eqref{eq:inner-control}.
\end{proof}

\begin{proof}[Proof of Theorem~\ref{thm:main}]
\emph{Necessity.} Formula \eqref{eq:solution} and the support properties of \eqref{eq:RL} give $S_Tq\in\AT$. The defect is linear in the state. Summing \eqref{eq:one-defect} over the inner sources gives \eqref{eq:cohom}.

\emph{Sufficiency.} Given the transfer functions $f_j$, use \eqref{eq:canonical-lift} and \eqref{eq:inner-control} to construct the inner controls. Subtract their states from $U$. The remainder is still in $\AT$ and its \emph{periodised} defect is zero by \eqref{eq:cohom}. Theorem~\ref{thm:two}, applied to $(z_1,z_M)$, then supplies $q_1,q_M$. This proves sufficiency and the reconstruction formula.
\end{proof}

\begin{remark}[Silent emissions]\label{rem:nonunique}
A nonzero control with zero final state will be called a \emph{silent emission}. The construction above identifies all possible sources of non-uniqueness. First, if there is an inner source and $cT>2d$, periodisation on $(0,cT)$ has a nontrivial kernel: for example, two opposite translates of the same function on disjoint intervals separated by $2d$ have zero periodisation. The state produced by such an inner control can be cancelled by the outer pair. Second, any nonzero solution of the homogeneous version of \eqref{eq:cohom} satisfying the support constraints can be lifted and cancelled in the same way. Conversely, if the final state is zero, its inner transfer functions solve this homogeneous equation; if these vanish and the lifts are unique, Theorem~\ref{thm:two} forces every control to vanish.
\end{remark}

\subsection{The defect quotient of the state space}
We make explicit the passage from statements about defects to statements about states. All circle norms below use Lebesgue measure on a period of length $2d$.

\begin{lemma}\label{lem:quotient}
The map $\GT:\YT\to L^2(\T_d)$ defined by \eqref{eq:defect} is bounded, and
\begin{equation}\label{eq:mean-defect}
\int_{\T_d}\GT U\dd y=\int_\R R\dd x-\int_\R L\dd x
=-c\int_\R u_0'\dd x.
\end{equation}
In particular, $\GT(\AT)\subset L^2_0(\T_d)$, the mean-zero subspace. If $cT\ge2d$, this map is onto $L^2_0(\T_d)$ and has a bounded right inverse $\JT$ satisfying
\begin{equation}\label{eq:J-norm}
\GT\JT g=g,\qquad
\norm{\JT g}_\X=\sqrt{\frac2k}\norm g_{L^2(\T_d)}.
\end{equation}
Its kernel on $\AT$ is the range of the outer pair. Thus, for $cT\ge2d$,
\begin{equation}\label{eq:quotient-decomp}
U=\big(U-\JT\GT U\big)+\JT\GT U,
\qquad U-\JT\GT U\in\ker\GT.
\end{equation}
For every linear subspace $D\subset L^2_0(\T_d)$,
\begin{equation}\label{eq:closure-preimage}
\overline{\GT^{-1}(D)}^{\AT}=\GT^{-1}(\overline D).
\end{equation}
\end{lemma}

\begin{proof}
The defect is supported in $[-d,d+cT]$. Periodisation has overlap multiplicity bounded, for example, by $2+m_T$. Since
\[
\norm{\Delta_U}_2^2\le2(\norm R_2^2+\norm L_2^2)=k\norm U_\X^2,
\]
one may use
\begin{equation}\label{eq:G-bound}
\norm{\GT U}_{L^2(\T_d)}\le C_G\norm U_\X,
\qquad C_G=\sqrt{k(2+m_T)}.
\end{equation}
Integrating the periodisation proves \eqref{eq:mean-defect}.

For $g\in L^2_0(\T_d)$, define
\[
F(x)=\ind_{(z_1,z_1+2d)}(x)g(x-z_1),\qquad
\JT g=\left(-\frac1c\int_{-\infty}^{x}F(t)\dd t,\ F(x)\right).
\]
The primitive vanishes on both sides of its support because $g$ has mean zero. The associated invariants are $R=F$ and $L=0$, so $\JT g\in\AT$ when $cT\ge2d$, and its periodised defect is $g$. Formula \eqref{eq:invariants} gives \eqref{eq:J-norm}. The assertion about the kernel follows from Theorem~\ref{thm:two}.

For \eqref{eq:closure-preimage}, continuity gives one inclusion. If $\GT U\in\overline D$, choose $g_m\in D$ converging to $\GT U$ and set
\[
U_m=U+\JT(g_m-\GT U).
\]
Then $\GT U_m=g_m$ and $U_m\to U$ in $\X$. This proves the reverse inclusion. In particular, density and closedness of a defect range transfer to its state preimage; nonclosedness also transfers, since $\JT$ maps convergent defect sequences to convergent state sequences.
\end{proof}

\section{Three sources: circle dynamics}\label{sec:three}

\subsection{The constrained rotation equation}
For $M=3$, equation \eqref{eq:cohom} has a single unknown $f=f_2$ and the rotation $y\mapsto y-\alpha$, where $\alpha=2d_2$. Its rotation number is $\theta=\alpha/(2d)=a_{12}/a_{13}$, and the equation is
\begin{equation}\label{eq:cohom3}
f(y-\alpha)-f(y)=\Pi(y),\qquad f=0\ \text{outside }A=A_2(T).
\end{equation}
For $T<2a_{13}$, let $J=\T_d\setminus A$, an arc of positive length $c(2a_{13}-T)$. Endpoints of arcs are immaterial throughout. For $y\in J$, define the first backward return time
\[
\nu(y)=\min\{n\ge1:y-n\alpha\in J\}.
\]
It is finite: irrational rotations meet every nonempty open arc, and rational rotations return after their period. An orbit disjoint from $J$ can occur only in the rational case. The geometry is illustrated in Figure~\ref{fig:circle}.

We fix the Fourier convention
\begin{equation}\label{eq:fourier-convention}
\widehat g_n=\frac1{2d}\int_0^{2d}g(y)e^{-i\pi ny/d}\dd y,
\qquad \norm g_2^2=2d\sum_{n\in\Z}|\widehat g_n|^2.
\end{equation}
For $s\ge0$, we use the equivalent Sobolev norm
\[
\norm g_{H^s(\T_d)}^2=2d\sum_{n\in\Z}\br n^{2s}|\widehat g_n|^2,
\qquad \br n=(1+n^2)^{1/2}.
\]
\begin{figure}[ht]
\centering
\begin{tikzpicture}[scale=1.15]
\def\r{1.6}
\draw[line width=6pt, blue!35] (200:\r) arc (200:470:\r);
\draw[line width=6pt, red!45] (110:\r) arc (110:200:\r);
\draw[thick] (0,0) circle (\r);
\node[blue!60!black,anchor=north west,align=left] at (-62:1.95) {arc $A=A_2(T)$\\[-1pt]\small($f$ may be nonzero)};
\node[red!70!black,anchor=south east,align=right] at (122:1.95) {arc $J$\\[-1pt]\small($f=0$)};
\foreach \k/\lab in {0/{y},1/{y-\alpha},2/{y-2\alpha},3/{y-3\alpha},4/{y-4\alpha},5/{y-5\alpha}} {
  \pgfmathsetmacro{\ang}{150-137.5*\k}
  \fill (\ang:\r) circle (1.7pt);
  \node at (\ang:\r+0.42) {\small $\lab$};
}
\foreach \k in {0,...,4} {
  \pgfmathsetmacro{\a}{150-137.5*\k}
  \pgfmathsetmacro{\b}{150-137.5*(\k+1)}
  \draw[->,gray!70,shorten <=3pt,shorten >=3pt] (\a:\r) -- (\b:\r);
}
\node[gray!50!black,font=\small,align=center] at (0,-0.15) {excursion\\ through $A$};
\end{tikzpicture}
\caption{The rotation of $\T_d$ by $-\alpha=-2(z_2-z_1)$. Starting from $y\in J$, the backward orbit makes an excursion through $A$ and returns to $J$ at $y-5\alpha$, so that $\nu(y)=5$. Here $\theta\approx0.382$ and $|J|=\frac14|\T_d|$, i.e.\ $T=\frac32a_{13}$. Reachability requires the sum of $\Pi$ along every such excursion to vanish.}
\label{fig:circle}
\end{figure}
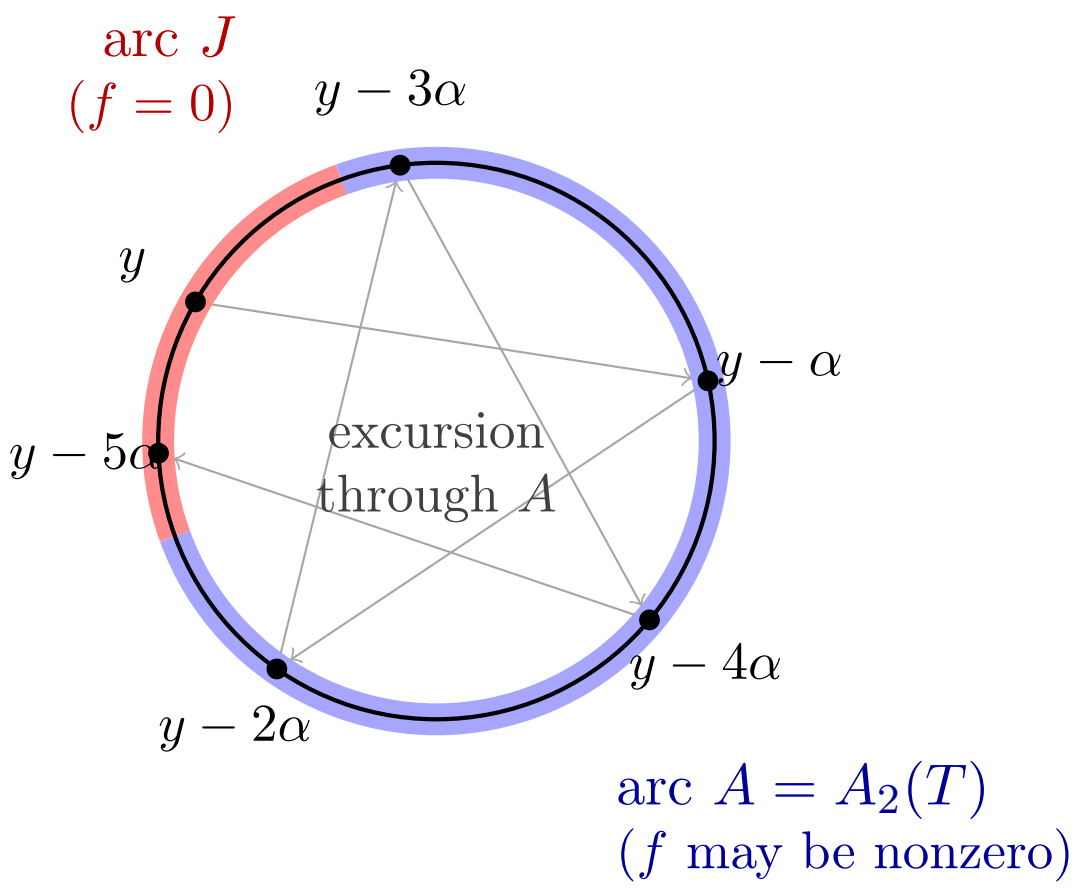
\begin{theorem}[Three sources]\label{thm:three}
Let $M=3$, let $U\in\AT$, and write $\Pi=\GT U$. When $\theta=m/v\in\Q$ is in lowest terms, put $\tau=a_{13}/v$. The state $U$ belongs to $\Rset$ if and only if the following corresponding conditions hold.
\begin{enumerate}[label=(\alph*)]
\item If $T<2a_{13}$, then
\begin{equation}\label{eq:excursion-condition}
\sum_{n=0}^{\nu(y)-1}\Pi(y-n\alpha)=0
\quad\text{for a.e. }y\in J.
\end{equation}
In the rational case, one additionally requires
\begin{equation}\label{eq:cycle-condition}
\sum_{n=0}^{v-1}\Pi(y-n\alpha)=0
\quad\text{for a.e. }y\text{ whose orbit is disjoint from }J.
\end{equation}
\item If $T\ge2a_{13}$, then
\begin{equation}\label{eq:three-fourier}
\widehat\Pi_n=0\quad(n\theta\in\Z),\qquad
\sum_{n\theta\notin\Z}\frac{|\widehat\Pi_n|^2}{4\sin^2(\pi n\theta)}<\infty.
\end{equation}
\end{enumerate}
Furthermore:
\begin{enumerate}[label=(\roman*)]
\item $\Rset$ is closed unless $\theta$ is irrational and $T\ge2a_{13}$. In that case $\overline{\Rset}=\AT$ and $\Rset\ne\AT$.
\item The controls are unique exactly when either $\theta$ is irrational and $T<2a_{13}$, or $\theta=m/v$ is rational and $T\le2(a_{13}-\tau)$.
\item If $\theta=m/v$ and $T\ge2a_{13}$, condition \eqref{eq:three-fourier} is equivalent to
\begin{equation}\label{eq:rational-periodisation}
\sum_{n\in\Z}\Delta_U(y+2nc\tau)=0
\quad\text{for a.e. }y.
\end{equation}
\end{enumerate}
\end{theorem}

\begin{proof}
\emph{The constrained equation.} Suppose $T<2a_{13}$. If $y\in A$ has a backward orbit meeting $J$, let
\[
N(y)=\min\{n\ge1:y-n\alpha\in J\}.
\]
Iterating \eqref{eq:cohom3} up to that first hit forces
\begin{equation}\label{eq:excursion-solution}
f(y)=-\sum_{n=0}^{N(y)-1}\Pi(y-n\alpha),\qquad f|_J=0.
\end{equation}
At a point of $A$ with $N(y)>1$, the same formula at $y-\alpha$ has one fewer term, and their difference is $\Pi(y)$. The conclusion is also immediate when $N(y)=1$. At $y\in J$, the equation becomes $f(y-\alpha)=\Pi(y)$, which is precisely \eqref{eq:excursion-condition}, including the case $\nu(y)=1$.

For irrational $\theta$, the translates of the interior of $J$ by the nonnegative rotation iterates cover the circle. Compactness gives a finite subcover, so all hitting and return times in this argument are bounded by a fixed $N_0=N_0(\theta,J)$. This is the elementary bounded-return consequence of minimality of an irrational rotation; see also \cite{KatokHasselblatt}. Formula \eqref{eq:excursion-solution} then defines a measurable function with
\[
\norm f_2\le N_0\norm\Pi_2.
\]

For $\theta=m/v$, every orbit is a cycle of length $v$. On the invariant measurable set $E$ of orbits disjoint from $J$, the equation is solvable exactly when \eqref{eq:cycle-condition} holds. A measurable particular solution on $E$ is
\begin{equation}\label{eq:cycle-solution}
f_E(y)=\frac1v\sum_{r=0}^{v-1}r\,\Pi(y-r\alpha),\qquad y\in E.
\end{equation}
Indeed, subtracting this expression from its translate yields $\Pi(y)$ minus the average of $\Pi$ over the cycle. On the complement of $E$, use \eqref{eq:excursion-solution}, whose sums have at most $v$ terms. The resulting function belongs to $L^2$. These constructions are bounded linear functions of $\Pi$ on the subspace defined by the finite orbit-sum constraints. Equivalently, each constraint is a bounded linear equation on the measurable sets where the relevant return time is fixed. There are only finitely many such sets. The admissible defect space is therefore closed, as is its preimage in $\AT$.

\emph{The unconstrained equation.} When $cT\ge2d$, there is no support constraint. Taking Fourier coefficients gives
\begin{equation}\label{eq:scalar-symbol}
\mu_n\widehat f_n=\widehat\Pi_n,
\qquad \mu_n=e^{-2\pi i n\theta}-1,
\qquad |\mu_n|=2|\sin(\pi n\theta)|.
\end{equation}
Parseval's identity gives exactly \eqref{eq:three-fourier}. For rational $\theta=m/v$, the nonzero multipliers are bounded below by $2\sin(\pi/v)$, so the range is the closed subspace with $\widehat\Pi_n=0$ for $n\in v\Z$.

For irrational $\theta$, only the constant mode is resonant. Every mean-zero trigonometric polynomial is in the range. Dirichlet's approximation principle gives nonzero integers for which $|\mu_n|$ is arbitrarily small. For completeness, choose distinct positive $n_r$ with $0<|\mu_{n_r}|\le2^{-r}$ and set
\[
\widehat g_{n_r}=\widehat g_{-n_r}=\frac{|\mu_{n_r}|}{\sqrt r},
\qquad \widehat g_n=0\ \text{otherwise}.
\]
Then $g\in L^2_0(\T_d)$, but the series in \eqref{eq:three-fourier} diverges. Its Fourier truncations are in the range, proving nonclosedness. Lemma~\ref{lem:quotient} transfers density and nonclosedness to $\AT$, proving (i).

\emph{Uniqueness.} For $cT<2d$, each lift is unique. A homogeneous solution of \eqref{eq:cohom3} is invariant under the rotation and vanishes on every orbit meeting $J$. For irrational $\theta$ all orbits meet $J$, so it vanishes. In the rational case, nonzero invariant functions supported on $E$ exist precisely when $E$ has positive measure. A rational orbit has $v$ equally spaced points, at distance $2d/v=2c\tau$. Thus $E$ has positive measure exactly when $|J|<2c\tau$, or $T>2(a_{13}-\tau)$. At equality, possible exceptional orbits form a null set. When $cT=2d$, a nonzero constant transfer function already gives a silent emission. For $cT>2d$, the same is true, and the lift also ceases to be unique. Remark~\ref{rem:nonunique} and uniqueness for the outer pair prove (ii).

\emph{The smaller period.} For rational $\theta=m/v$, the resonant Fourier indices are exactly $v\Z$. If $h=2d/v=2c\tau$, the $h$-periodisation of $\Delta_U$ has Fourier coefficients
\[
\frac1h\int_\R\Delta_U(y)e^{-2\pi i r y/h}\dd y
   =v\widehat\Pi_{rv},\qquad r\in\Z.
\]
Its vanishing is therefore equivalent to the resonant conditions, proving (iii).
\end{proof}

\begin{example}\label{ex:silent}
For any three-source geometry and $T=2a_{13}$, the controls
\[
\begin{aligned}
q_1&=\ind_{[a_{12},\,a_{12}+2a_{23}]},\\
q_2&=-\ind_{[0,\,2a_{13}]},\\
q_3&=\ind_{[a_{23},\,a_{23}+2a_{12}]}
\end{aligned}
\]
produce the zero final state. Substitution in \eqref{eq:rl} gives both $r=0$ and $l=0$; the middle transfer function is the constant $-\rho c/2$. For equally spaced sources, with $a=a_{12}=a_{23}$, silent emissions exist already for $T>2a$: choose a nonzero $\psi\in L^2(0,T-2a)$, extend it by zero, and set
\[
(q_1,q_2,q_3)=\big(\psi(\cdot-a),\,-\psi-\psi(\cdot-2a),\,\psi(\cdot-a)\big).
\]
All three controls are supported in $[0,T]$, and their translated contributions cancel in \eqref{eq:rl}.
\end{example}

\subsection{Diophantine regularity and spectral cost}
Let $\dist{x}=\min_{m\in\Z}|x-m|$. We will use
\begin{equation}\label{eq:sine-distance}
2\dist{x}\le |\sin(\pi x)|\le\pi\dist{x}.
\end{equation}
For $T\ge2a_{13}$ and a defect satisfying the resonant conditions, define
\begin{equation}\label{eq:Dtheta}
\DK_\theta(g)=\left(2d\sum_{n\theta\notin\Z}
\frac{|\widehat g_n|^2}{4\sin^2(\pi n\theta)}\right)^{1/2}.
\end{equation}
By \eqref{eq:scalar-symbol}, this is exactly the minimum $L^2(\T_d)$ norm of a transfer function solving the circle equation. The minimiser has zero coefficients on the resonant modes.

\begin{corollary}\label{cor:dioph}
Let $\theta$ be irrational and $T\ge2a_{13}$.
\begin{enumerate}[label=(\roman*)]
\item If $\dist{n\theta}\ge b|n|^{-\nu}$ for all $n\ne0$, then every $U\in\AT$ with $\Pi_U\in H^\nu(\T_d)$ is reachable, and
\begin{equation}\label{eq:scalar-regularity}
\DK_\theta(\Pi_U)\le\frac1{4b}\norm{\Pi_U}_{H^\nu(\T_d)}.
\end{equation}
In particular, $\nu=1$ suffices for badly approximable numbers, including quadratic irrationals.
\item If $\theta$ is Liouville, there are states in $\AT$ with $C^\infty$ periodised defect which are not reachable.
\end{enumerate}
\end{corollary}

\begin{proof}
The first statement follows from \eqref{eq:sine-distance} and \eqref{eq:Dtheta}. The stated facts about badly approximable and quadratic irrational numbers are standard results of continued-fraction theory; see \cite{Cassels}.

For the second, choose increasing positive integers $n_r\ge2^r$ for which
$\delta_r:=\dist{n_r\theta}<n_r^{-r}$. Define a real mean-zero periodic function by
\[
\widehat g_{n_r}=\widehat g_{-n_r}=\delta_r^{1/2},\qquad
\widehat g_n=0\ \text{otherwise}.
\]
For each $s\ge0$, the tail of $\sum_r\br{n_r}^{2s}\delta_r$ is bounded by a constant times $\sum_r n_r^{-2}$, hence is finite. Thus $g\in H^s$ for every $s$, and $g\in C^\infty$. By the upper bound in \eqref{eq:sine-distance}, its terms in \eqref{eq:Dtheta} are bounded below by a positive constant times $\delta_r^{-1}$, so $\DK_\theta(g)=\infty$. The state $U=\JT g$ has the required defect and is not reachable. Smoothness is asserted here for the periodic defect, not for its zero-extended state representative at the cut endpoints.
\end{proof}

\begin{proposition}\label{prop:scalar-cost}
For $N\ge1$, set
\begin{equation}\label{eq:kappa-scalar}
\kappa_N(\theta)=
\max_{\substack{1\le |n|\le N\\n\theta\notin\Z}}
\frac1{2|\sin(\pi n\theta)|}.
\end{equation}
On the subspace of trigonometric polynomials of degree at most $N$ with vanishing resonant coefficients, $\kappa_N(\theta)$ is the exact norm of the least-norm inverse of \eqref{eq:cohom3}. More precisely,
\[
\sup_{g\ne0}\frac{\DK_\theta(g)}{\norm g_2}=\kappa_N(\theta)
\]
on this subspace. The following arithmetic estimates hold.
\begin{enumerate}[label=(\roman*)]
\item If $\theta=m/v$ in lowest terms, then $\kappa_N(\theta)\le[2\sin(\pi/v)]^{-1}$, with equality for $N\ge v-1$.
\item For every irrational $\theta$, $\kappa_N(\theta)\ge N/(2\pi)$. If also $\dist{n\theta}\ge b/|n|$, then
\begin{equation}\label{eq:badly-cost}
\frac{N}{2\pi}\le\kappa_N(\theta)\le\frac{N}{4b}
\qquad(N\ge1).
\end{equation}
More generally, a lower bound $\dist{n\theta}\ge b|n|^{-\nu}$ gives $\kappa_N(\theta)\le(4b)^{-1}N^\nu$.
\item If $\theta$ is Liouville, then for every $r>0$,
\begin{equation}\label{eq:liouville-growth}
\limsup_{N\to\infty}\frac{\kappa_N(\theta)}{N^r}=\infty.
\end{equation}
\end{enumerate}
\end{proposition}

\begin{proof}
The norm identity follows from the diagonal Fourier formula; equality is achieved by a real sine or cosine at a maximising nonresonant frequency. For rational $\theta=m/v$, multiplication by $m$ permutes the nonzero residues modulo $v$, so the least nonzero sine is $\sin(\pi/v)$ and is attained by an integer between $1$ and $v-1$.

For irrational $\theta$, the pigeonhole principle applied to $0,\theta,\dots,N\theta$ modulo one gives an integer $1\le n\le N$ with $\dist{n\theta}\le1/N$. Using the upper sine bound proves the lower estimate in \eqref{eq:badly-cost}. The lower sine bound gives the upper estimates. Finally, for a Liouville number there are arbitrarily large $n$ with $\dist{n\theta}<n^{-m}$, for any prescribed $m$. Taking $N=n$ gives $\kappa_N\ge(2\pi)^{-1}N^m$ along such choices, which proves \eqref{eq:liouville-growth}.
\end{proof}

\begin{example}\label{ex:arithmetic}
Fix $z_1,z_3$ and take $z_2=z_1+\theta(z_3-z_1)$.
For $\theta=1/2$, exactly the even Fourier modes are obstructed, and $\kappa_N=1/2$ on the accessible modes.
For $\theta=(\sqrt5-1)/2$, all nonzero modes are accessible, and \eqref{eq:badly-cost} gives linear growth. One can use $b=1/3$: if $p$ is the nearest integer to $n\theta$, then
\[
1\le|p^2+pn-n^2|
 =n^2|p/n-\theta|\,|p/n+(1+\sqrt5)/2|
 <3n\,|p-n\theta|.
\]
For the Liouville number $\theta=\sum_{r=1}^\infty10^{-r!}$, let $N_m=10^{m!}$, $m\ge2$. Truncating the series gives
\[
\dist{N_m\theta}<2N_m^{-m},\qquad
\kappa_{N_m}(\theta)>\frac{N_m^m}{4\pi}.
\]
Thus rational placement creates exact inaccessible modes, a badly approximable placement avoids such modes with linear cut-off cost, and Liouville placement avoids exact resonances but can require much larger controls.
\end{example}

The Fourier structure here is analogous to arithmetic effects for bounded strings, but its origin is the cancellation of outgoing states on the line. The next section proves that these inverse norms control the actual actuator norm and gives the corresponding approximation estimates.

\section{Long times for \texorpdfstring{$M$}{M} sources}\label{sec:long}

Throughout this section, $M\ge3$ and $T\ge2a_{1M}$, so every arc $A_j(T)$ is the whole circle. Write $p=M-2$, $\boldsymbol\theta=(\theta_2,\dots,\theta_{M-1})$, and define
\begin{equation}\label{eq:multi-symbol}
\mu_{j,n}=e^{-2\pi i n\theta_j}-1,\qquad
\sigma_n=\left(\sum_{j=2}^{M-1}|\mu_{j,n}|^2\right)^{1/2}
        =\left(4\sum_{j=2}^{M-1}\sin^2(\pi n\theta_j)\right)^{1/2}.
\end{equation}
The common resonant set and the corresponding closed defect space are
\begin{equation}\label{eq:resonant-space}
\begin{aligned}
\Kres&=\{n\in\Z:n\theta_j\in\Z\text{ for every }1<j<M\},\\
\mathcal H_{\rm res}^{\perp}
&=\{g\in L^2(\T_d):\widehat g_n=0\text{ for }n\in\Kres\}.
\end{aligned}
\end{equation}
In particular, $0\in\Kres$ and this defect space consists of mean-zero functions.

\subsection{The long-time reachable set}
\begin{theorem}[Long times]\label{thm:long}
Let $U\in\AT$ and $g=\GT U$. Then $U\in\Rset$ if and only if $g\in\mathcal H_{\rm res}^{\perp}$ and
\begin{equation}\label{eq:multi-D}
\DK_{\boldsymbol\theta}(g)^2
 :=2d\sum_{n\notin\Kres}\frac{|\widehat g_n|^2}{\sigma_n^2}<\infty.
\end{equation}
The transfer functions of least total squared norm are given by
\begin{equation}\label{eq:multi-minimizer}
\widehat f_{j,n}^{\,\min}=
\begin{cases}
\displaystyle\frac{\overline{\mu_{j,n}}}{\sigma_n^2}\widehat g_n,
       &n\notin\Kres,\\[4pt]
0,&n\in\Kres,
\end{cases}
\end{equation}
and satisfy
\begin{equation}\label{eq:transfer-cost}
\sum_{j=2}^{M-1}\norm{f_j^{\min}}_2^2
 =\DK_{\boldsymbol\theta}(g)^2.
\end{equation}
Moreover,
\begin{enumerate}[label=(\roman*)]
\item The reachable set is closed if and only if all $\theta_j$ are rational. In that case, if $V$ is their least common denominator and $\tau=a_{1M}/V$, then $\Kres=V\Z$ and reachability is equivalent to
\begin{equation}\label{eq:multi-rational}
\sum_{n\in\Z}\Delta_U(y+2nc\tau)=0
\quad\text{for a.e. }y.
\end{equation}
\item In general,
\begin{equation}\label{eq:closure-long}
\overline{\Rset}=\{U\in\AT:\GT U\in\mathcal H_{\rm res}^{\perp}\}.
\end{equation}
Consequently, $\Rset$ is dense in $\AT$ if and only if at least one $\theta_j$ is irrational. In that case it is not closed.
\end{enumerate}
\end{theorem}

\begin{proof}
In Fourier variables, \eqref{eq:cohom} is the scalar equation
\[
\widehat g_n=\sum_{j=2}^{M-1}\mu_{j,n}\widehat f_{j,n}.
\]
If $n\in\Kres$, it is solvable exactly when $\widehat g_n=0$. Otherwise Cauchy--Schwarz gives
\[
\sum_j|\widehat f_{j,n}|^2\ge\frac{|\widehat g_n|^2}{\sigma_n^2},
\]
with equality for \eqref{eq:multi-minimizer}. Parseval's identity proves necessity and sufficiency of \eqref{eq:multi-D}, as well as \eqref{eq:transfer-cost}. For real $g$, the coefficients in \eqref{eq:multi-minimizer} have the required conjugate symmetry, so the transfer functions and the controls can be chosen real.

If all the angles are rational, $\sigma_n$ is periodic in $n$ and has a positive minimum away from its zeros. The common zero set is $V\Z$. The same Fourier calculation used for \eqref{eq:rational-periodisation} gives \eqref{eq:multi-rational}. Thus the defect range equals $\mathcal H_{\rm res}^{\perp}$ and is closed.

If some angle is irrational, $\Kres=\{0\}$. To see that $\inf_{n\ne0}\sigma_n=0$, divide $[0,1)^p$ into $Q^p$ cubes of side $1/Q$ and apply the pigeonhole principle to the $Q^p+1$ points $r\boldsymbol\theta$ modulo $\Z^p$, $0\le r\le Q^p$. Their difference gives an integer $1\le n\le Q^p$ with
\begin{equation}\label{eq:dirichlet-multi}
\max_j\dist{n\theta_j}\le Q^{-1},\qquad
\sigma_n\le2\pi\sqrt p\,Q^{-1}.
\end{equation}
Since there is no nonzero common resonant index, such integers are unbounded as $Q\to\infty$. Choosing a subsequence with $\sigma_{n_r}\le2^{-r}$ and coefficients $\widehat g_{\pm n_r}=\sigma_{n_r}/\sqrt r$ gives an $L^2$ defect not satisfying \eqref{eq:multi-D}. On the other hand, every finite Fourier sum with vanishing resonant coefficients does satisfy it. Hence the defect range is dense in $\mathcal H_{\rm res}^{\perp}$ and is nonclosed in the irrational case. Lemma~\ref{lem:quotient} yields \eqref{eq:closure-long} and all remaining assertions. The rational defect space is proper because it excludes nonzero modes in $V\Z$.
\end{proof}

Once $T\ge2a_{1M}$, the criterion on the periodised defect is independent of $T$. The admissible state space and the constants in the control estimates below still depend on $T$.

\subsection{From transfer functions to the physical control cost}
We extend $\CT(U)$ by $+\infty$ for unreachable states. For reachable states its minimum exists: $S_T^{-1}\{U\}$ is a nonempty closed affine subspace of the Hilbert space $L^2(0,T)^M$. The minimiser is unique, even when there are silent emissions, because it is the unique element orthogonal to $\ker S_T$ in that affine subspace.

\begin{theorem}[Two-sided control-cost comparison]\label{thm:physical-cost}
For fixed $T,\rho,c$ and source positions, there are constants $0<c_T\le C_T<\infty$ such that every $U\in\Rset$, with $g=\GT U$, satisfies
\begin{equation}\label{eq:physical-equivalence}
c_T\big(\norm U_\X+\DK_{\boldsymbol\theta}(g)\big)
\le\CT(U)
\le C_T\big(\norm U_\X+\DK_{\boldsymbol\theta}(g)\big).
\end{equation}
More explicitly, with $\beta,m_T$ from \eqref{eq:lift-constants} and $B_T$ from \eqref{eq:outer-bound} for the outer pair,
\begin{equation}\label{eq:physical-lower}
\CT(U)\ge
\max\left\{\frac{\norm U_\X}{\sqrt{M\rho c}},\,
\frac{\beta}{\sqrt{m_T}}\DK_{\boldsymbol\theta}(g)\right\},
\end{equation}
and a control obtained from \eqref{eq:multi-minimizer}, the single-period lift, and the outer recursion satisfies
\begin{equation}\label{eq:physical-upper}
\norm q_{L^2(0,T)^M}\le
B_T\norm U_\X+
\beta\big(1+B_T\sqrt{p\rho c}\big)\DK_{\boldsymbol\theta}(g).
\end{equation}
\end{theorem}

\begin{proof}
For any control reaching $U$, Proposition~\ref{prop:one} gives the first lower bound. Its inner transfer functions solve the circle equation and therefore have norm at least $\DK_{\boldsymbol\theta}(g)$. By \eqref{eq:period-bound},
\[
\DK_{\boldsymbol\theta}(g)^2
\le\sum_{j=2}^{M-1}\norm{f_j}_2^2
\le\frac{m_T}{\beta^2}\sum_{j=2}^{M-1}\norm{q_j}_2^2.
\]
This proves the second lower bound for every admissible control.

For the upper bound choose the least-norm transfer functions \eqref{eq:multi-minimizer}. Their single-period lifts give inner controls $q^{\rm in}$ with
\[
\norm{q^{\rm in}}_2=\beta\DK_{\boldsymbol\theta}(g).
\]
Let $W$ be their total final state. The single-source estimate and Cauchy--Schwarz over the $p$ inner sources imply
\[
\norm W_\X\le\sqrt{p\rho c}\,\beta\DK_{\boldsymbol\theta}(g).
\]
The remainder $U-W$ has zero periodised defect. The outer recursion gives $q^{\rm out}$ satisfying
\[
\norm{q^{\rm out}}_2\le B_T\norm{U-W}_\X.
\]
Adding the two control norms gives \eqref{eq:physical-upper}. Finally, the maximum of the two positive quantities in \eqref{eq:physical-lower} bounds a fixed positive multiple of their unweighted sum, proving \eqref{eq:physical-equivalence}.
\end{proof}

\begin{remark}\label{rem:cost-meaning}
The norm \eqref{eq:transfer-cost} is an exact minimum for the auxiliary circle problem. Formula \eqref{eq:physical-upper} constructs a physical control with comparable norm; it is not, in general, a formula for the least-norm physical control, because changing the inner transfer functions can also change the cost of the outer correction. The distinction is immaterial for the two-sided estimates but is essential in interpreting them.
\end{remark}

\subsection{Spectral synthesis and approximation of a prescribed state}
Let $P_N$ be the Fourier projection onto $|n|\le N$, and put
\begin{equation}\label{eq:multi-kappa}
\begin{aligned}
\mathcal V_N&=\{g\in\mathcal H_{\rm res}^{\perp}:P_Ng=g\},\\
\kappa_N(\boldsymbol\theta)&=
\max_{\substack{1\le |n|\le N\\n\notin\Kres}}\sigma_n^{-1}.
\end{aligned}
\end{equation}
The maximum is over a nonempty set for $N\ge1$, since each inner angle lies strictly between zero and one.

\begin{theorem}[Band-limited cost and constructive approximation]\label{thm:cutoff}
The following statements hold.
\begin{enumerate}[label=(\roman*)]
\item The exact inverse norm for defects in $\mathcal V_N$ is
\begin{equation}\label{eq:inverse-norm}
\sup_{0\ne g\in\mathcal V_N}
\frac{\DK_{\boldsymbol\theta}(g)}{\norm g_2}
=\kappa_N(\boldsymbol\theta).
\end{equation}
Define the worst physical synthesis cost on states with defects of degree at most $N$ by
\[
\mathfrak K_T(N)=
\sup_{\substack{0\ne U\in\AT\\\GT U\in\mathcal V_N}}
\frac{\CT(U)}{\norm U_\X}.
\]
There are constants $a_T,A_T>0$, independent of $N$, such that
\begin{equation}\label{eq:physical-cutoff}
a_T\big(1+\kappa_N(\boldsymbol\theta)\big)
\le\mathfrak K_T(N)
\le A_T\big(1+\kappa_N(\boldsymbol\theta)\big).
\end{equation}
\item Let $U\in\AT$ and $g=\GT U\in\mathcal H_{\rm res}^{\perp}$. Define
\begin{equation}\label{eq:UN}
U_N=U-\JT(g-P_Ng).
\end{equation}
Then $U_N\in\Rset$, $\GT U_N=P_Ng$, and there are explicitly reconstructible controls $q^{(N)}$ with $S_Tq^{(N)}=U_N$ such that
\begin{align}
\norm{U-U_N}_\X
 &=\sqrt{\frac2k}\norm{(\id-P_N)g}_2,
\label{eq:approx-error}\\
\norm{q^{(N)}}_2
 &\le C_T\big(\norm U_\X+\DK_{\boldsymbol\theta}(P_Ng)\big)
\nonumber\\
 &\le C_T\big(\norm U_\X+\kappa_N(\boldsymbol\theta)\norm g_2\big).
\label{eq:approx-cost}
\end{align}
In particular, $U_N\to U$ in energy. If $g\in H^s(\T_d)$, $s>0$, then
\begin{equation}\label{eq:sobolev-error}
\norm{U-U_N}_\X\le\sqrt{\frac2k}\br{N+1}^{-s}\norm g_{H^s}.
\end{equation}
\item Suppose, on all nonresonant indices,
\begin{equation}\label{eq:dioph-sigma}
\sigma_n\ge\gamma|n|^{-\nu},\qquad \gamma>0,\quad\nu\ge0.
\end{equation}
For $g\in H^s$, $s\ge0$, the controls above satisfy the sharper bound
\begin{equation}\label{eq:sobolev-cost}
\norm{q^{(N)}}_2\le C_T\left(\norm U_\X+
\gamma^{-1}N^{(\nu-s)_+}\norm g_{H^s}\right).
\end{equation}
If $s\ge\nu$, then $U$ itself is reachable with
\begin{equation}\label{eq:exact-sobolev-cost}
\CT(U)\le C_T\left(\norm U_\X+\gamma^{-1}\norm g_{H^s}\right).
\end{equation}
If $0<s<\nu$, $g\ne0$, and $0<\varepsilon\le\norm g_{H^s}$, there is a control $q_\varepsilon$ such that
\begin{equation}\label{eq:epsilon-cost}
\begin{aligned}
\norm{S_Tq_\varepsilon-U}_\X&\le\varepsilon,\\
\norm{q_\varepsilon}_2&\le C\left(\norm U_\X+
\norm g_{H^s}^{\nu/s}\varepsilon^{-(\nu-s)/s}\right),
\end{aligned}
\end{equation}
where $C$ may also depend on $s,\nu,\gamma$.
\end{enumerate}
\end{theorem}

\begin{proof}
The diagonal formula \eqref{eq:multi-D} proves \eqref{eq:inverse-norm}; a real normalised cosine at a maximising frequency attains equality. By Theorem~\ref{thm:physical-cost} and \eqref{eq:G-bound},
\[
\CT(U)\le C_T(1+C_G\kappa_N)\norm U_\X
\]
when $\GT U\in\mathcal V_N$. For the reverse inequality, take such a normalised cosine $g$ and $U=\JT g$. Then \eqref{eq:physical-lower} and \eqref{eq:J-norm} give
\[
\frac{\CT(U)}{\norm U_\X}
\ge\frac{\beta\sqrt{k/2}}{\sqrt{m_T}}\kappa_N.
\]
The first lower bound in \eqref{eq:physical-lower} also gives a positive $N$-independent bound for $\mathfrak K_T(N)$; the admissible set is nonempty, for instance because it contains nonzero outer-source states. Combining the two bounds proves \eqref{eq:physical-cutoff}.

For (ii), $\GT\JT=\id$ gives $\GT U_N=P_Ng$. This is a finite Fourier sum with no resonant component, so Theorem~\ref{thm:long} implies reachability. Formula \eqref{eq:J-norm} gives the error identity. Use the controls constructed in Theorem~\ref{thm:physical-cost} for the target $U_N$. Since
\[
\norm{U_N}_\X\le\norm U_\X+\sqrt{2/k}\norm g_2
\le(1+\sqrt{2/k}C_G)\norm U_\X,
\]
we obtain \eqref{eq:approx-cost}. The Sobolev error bound follows by estimating the Fourier tail.

For (iii), \eqref{eq:dioph-sigma} gives
\[
\begin{split}
\DK_{\boldsymbol\theta}(P_Ng)^2
&\le\frac{2d}{\gamma^2}
 \sum_{\substack{1\le|n|\le N\\n\notin\Kres}}
 |n|^{2\nu}|\widehat g_n|^2\\
&\le\gamma^{-2}N^{2(\nu-s)_+}\norm g_{H^s}^2.
\end{split}
\]
This proves \eqref{eq:sobolev-cost}. When $s\ge\nu$, the same estimate on the full series proves exact reachability and \eqref{eq:exact-sobolev-cost}. When $0<s<\nu$, choose an integer $N$ comparable to $(\norm g_{H^s}/\varepsilon)^{1/s}$, with the constant large enough for \eqref{eq:sobolev-error}. Substituting this choice into \eqref{eq:sobolev-cost} gives \eqref{eq:epsilon-cost}.
\end{proof}

\begin{remark}
The cut-off is imposed on the \emph{periodised defect}, not on the Fourier transform of the entire spatial state. The part $U-\JT\GT U$ is retained exactly and handled by the outer pair. Thus \eqref{eq:UN} approximates an arbitrary target in the closure of the reachable set without discarding its zero-defect component.
\end{remark}

\subsection{Simultaneous Diophantine approximation and actuator geometry}
\begin{proposition}\label{prop:multi-arithmetic}
Suppose at least one $\theta_j$ is irrational. Then for every $N\ge1$,
\begin{equation}\label{eq:universal-multi-lower}
\kappa_N(\boldsymbol\theta)\ge\frac{N^{1/p}}{4\pi\sqrt p}.
\end{equation}
If the vector is simultaneously badly approximable in the sense that
\begin{equation}\label{eq:simultaneous-BA}
\max_{2\le j\le M-1}\dist{n\theta_j}
\ge b|n|^{-1/p}\qquad(n\in\Z\setminus\{0\}),
\end{equation}
then
\begin{equation}\label{eq:multi-sharp}
\frac{N^{1/p}}{4\pi\sqrt p}
\le\kappa_N(\boldsymbol\theta)
\le\frac{N^{1/p}}{4b}.
\end{equation}
Conversely, a bound $\kappa_N(\boldsymbol\theta)\le C N^{1/p}$ for all $N$ implies \eqref{eq:simultaneous-BA} with a possibly different constant. Consequently, under \eqref{eq:simultaneous-BA},
\[
\mathfrak K_T(N)\asymp N^{1/p}.
\]
\end{proposition}

\begin{proof}
Choose $Q=\lfloor N^{1/p}\rfloor$. Then $Q\ge N^{1/p}/2$, and \eqref{eq:dirichlet-multi} supplies $1\le n\le Q^p\le N$ with
$\sigma_n\le4\pi\sqrt p\,N^{-1/p}$. It is nonresonant because some angle is irrational. This proves \eqref{eq:universal-multi-lower}. By \eqref{eq:sine-distance},
\[
4\max_j\dist{n\theta_j}\le\sigma_n
\le2\pi\sqrt p\max_j\dist{n\theta_j}.
\]
The first inequality gives the upper bound in \eqref{eq:multi-sharp}. Conversely, set $N=|n|$ in the assumed estimate: $\sigma_n^{-1}\le C|n|^{1/p}$. The second inequality then yields \eqref{eq:simultaneous-BA}. The physical cost assertion follows from \eqref{eq:physical-cutoff}.
\end{proof}

\begin{remark}[More sources do not automatically improve the exponent]\label{rem:dependent-angles}
For four sources, take $\theta_2=\alpha$, $\theta_3=2\alpha$, where $0<\alpha<1/2$ is badly approximable. Since $|\sin(2x)|\le2|\sin x|$,
\[
2|\sin(\pi n\alpha)|\le\sigma_n
\le2\sqrt5\,|\sin(\pi n\alpha)|.
\]
Thus $\kappa_N(\alpha,2\alpha)\asymp N$, rather than $N^{1/2}$. The rational dependence of the two angles preserves a one-dimensional small-divisor obstruction. In contrast, \eqref{eq:simultaneous-BA} is precisely the condition needed for the rate $N^{1/p}$.
\end{remark}

\subsection{Sharp worst-case cost at prescribed accuracy}
To distinguish an upper bound for one reconstruction algorithm from an unavoidable cost, we now allow \emph{all} controls in the lower estimate. Put
\begin{equation}\label{eq:approx-min-cost}
\mathfrak C_T^\varepsilon(U)
 =\inf\{\norm q_{L^2(0,T)^M}:\norm{S_Tq-U}_\X\le\varepsilon\}.
\end{equation}
For $s>0$ and $\Kres=\{0\}$, define
\begin{equation}\label{eq:minimax-def}
\mathfrak E_s(\varepsilon)
 =\sup_{\substack{g\in H^s(\T_d),\ \widehat g_0=0\\\norm g_{H^s}\le1}}
 \mathfrak C_T^\varepsilon(\JT g).
\end{equation}
This is the worst approximation cost on the canonical state representatives of the unit Sobolev ball of defects. Restricting to these representatives removes the freely specifiable zero-defect component. Adding any uniformly energy-bounded zero-defect component changes the upper estimate only by an additive bounded amount, by Theorem~\ref{thm:two}.

\begin{lemma}\label{lem:mode-lower}
Let $n\ge1$ be nonresonant and let
\[
h_n(y)=d^{-1/2}\cos(\pi ny/d),\qquad \norm{h_n}_2=1.
\]
For $a>0$ and any control satisfying
$\norm{S_Tq-\JT(ah_n)}_\X\le\varepsilon$, one has
\begin{equation}\label{eq:mode-lower}
\norm q_2\ge\frac{\beta}{\sqrt{m_T}}\,
\frac{(a-C_G\varepsilon)_+}{\sigma_n},
\end{equation}
where $C_G$ is as in \eqref{eq:G-bound}.
\end{lemma}

\begin{proof}
Let $g_q=\GT S_Tq$. The boundedness of $\GT$ gives
$\norm{g_q-ah_n}_2\le C_G\varepsilon$. Orthogonal projection onto the two modes $\{-n,n\}$ therefore has norm at least $(a-C_G\varepsilon)_+$. Since $\sigma_{-n}=\sigma_n$, \eqref{eq:multi-D} yields
\[
\DK_{\boldsymbol\theta}(g_q)
\ge\frac{(a-C_G\varepsilon)_+}{\sigma_n}.
\]
The second bound in \eqref{eq:physical-lower}, applied to the reachable state $S_Tq$, proves the result. No assumption is made on the other Fourier modes of the state produced by $q$.
\end{proof}

\begin{theorem}[Sharp Sobolev cost--accuracy rates]\label{thm:minimax}
Assume \eqref{eq:simultaneous-BA}, put $\nu=1/p$, and let $0<s<\nu$. There are constants $c,C,\varepsilon_0>0$ such that
\begin{equation}\label{eq:minimax-sharp}
c\,\varepsilon^{-(\nu-s)/s}
\le\mathfrak E_s(\varepsilon)
\le C\,\varepsilon^{-(\nu-s)/s}
\qquad(0<\varepsilon<\varepsilon_0).
\end{equation}
The lower bound in \eqref{eq:minimax-sharp} remains valid whenever at least one angle is irrational, without the lower Diophantine bound \eqref{eq:simultaneous-BA}. For each $0<s<\nu$, there are mean-zero defects in $H^s(\T_d)$ whose canonical states are not exactly reachable. For $s\ge\nu$ under \eqref{eq:simultaneous-BA}, the canonical Sobolev class is exactly reachable with uniformly bounded control cost.

In particular, for three sources with badly approximable irrational $\theta$, \eqref{eq:minimax-sharp} becomes
\begin{equation}\label{eq:scalar-minimax}
\mathfrak E_s(\varepsilon)\asymp\varepsilon^{-(1-s)/s},
\qquad 0<s<1.
\end{equation}
For three sources with Liouville $\theta$, for every $s>0$ and $r>0$,
\begin{equation}\label{eq:liouville-minimax}
\limsup_{\varepsilon\downarrow0}
\varepsilon^r\mathfrak E_s(\varepsilon)=\infty.
\end{equation}
\end{theorem}

\begin{proof}
Under \eqref{eq:simultaneous-BA}, $\sigma_n\ge4b|n|^{-\nu}$. For $\norm g_{H^s}\le1$, the state $\JT g$ has uniformly bounded energy norm. Choose $N$ comparable to $\varepsilon^{-1/s}$ in \eqref{eq:sobolev-error} and \eqref{eq:sobolev-cost}. The error is at most $\varepsilon$ and the control norm is at most $C\varepsilon^{-(\nu-s)/s}$ for small $\varepsilon$. Taking the supremum proves the upper bound. The assertion for $s\ge\nu$ follows from \eqref{eq:exact-sobolev-cost}.

For the lower bound, fix small $\varepsilon$ and choose an integer $N\asymp\varepsilon^{-1/s}$ such that
\[
\br N^{-s}\ge2C_G\varepsilon.
\]
The pigeonhole argument in Proposition~\ref{prop:multi-arithmetic} supplies $1\le n\le N$ with $0<\sigma_n\le C_pN^{-\nu}$. Set
\[
a=\br n^{-s},\qquad g=ah_n.
\]
Then $\norm g_{H^s}=1$, and $a\ge\br N^{-s}\ge2C_G\varepsilon$. Lemma~\ref{lem:mode-lower} applies to every control approximating $\JT g$ to error $\varepsilon$ and gives
\[
\norm q_2\ge\frac{\beta}{2\sqrt{m_T}}\frac{a}{\sigma_n}
\ge cN^{\nu-s}
\ge c'\varepsilon^{-(\nu-s)/s}.
\]
This proof uses only the absence of nonzero common resonances, not \eqref{eq:simultaneous-BA}.

To verify sharpness of the exact regularity threshold, choose distinct positive integers $n_j\to\infty$ for which
\[
\sigma_{n_j}\le C_p n_j^{-\nu},\qquad n_j^{\nu-s}\ge j^2.
\]
Such a subsequence follows from \eqref{eq:dirichlet-multi}: the indices are unbounded, and $n\le Q^p$ implies $Q^{-1}\le n^{-1/p}$. Set
\[
g=\sum_{j\ge1}j^{-1}\br{n_j}^{-s}h_{n_j}.
\]
Orthogonality gives $\norm g_{H^s}^2=\sum_j j^{-2}<\infty$, whereas
\[
\DK_{\boldsymbol\theta}(g)^2
 =\sum_{j\ge1}\frac{j^{-2}\br{n_j}^{-2s}}{\sigma_{n_j}^2}
 \ge c\sum_{j\ge1}j^{-2}n_j^{2(\nu-s)}=\infty.
\]
Thus $\JT g$ is not exactly reachable.

Finally, let $p=1$ and $\theta$ be Liouville. Choose increasing $n_j\to\infty$ with $\dist{n_j\theta}<n_j^{-j}$. Set
\[
a_j=\br{n_j}^{-s},\qquad
\varepsilon_j=\frac{a_j}{2C_G},\qquad g_j=a_jh_{n_j}.
\]
Then $\norm{g_j}_{H^s}=1$, and \eqref{eq:mode-lower}, together with $\sigma_{n_j}\le2\pi n_j^{-j}$, yields
\[
\mathfrak E_s(\varepsilon_j)\ge c\,n_j^{j-s}.
\]
For every fixed $r>0$,
\[
\varepsilon_j^r\mathfrak E_s(\varepsilon_j)
\ge c_r n_j^{j-s(r+1)}\longrightarrow\infty,
\]
which proves \eqref{eq:liouville-minimax}.
\end{proof}

These estimates distinguish three effects of actuator geometry. Rational relations create exact inaccessible Fourier modes. Badly approximable placements remove those exact obstructions while retaining a controlled algebraic amplification, with a sharp regularity threshold and approximation cost. Liouville placements preserve density but allow arbitrarily severe loss of stability. The outer-pair reconstruction shows that these are properties of the original point-source control problem, rather than artefacts of its circle representation.
\bigskip

\noindent\textbf{AI Declaration}\par

\noindent AI tools were used to assist with language editing, mathematical review, and the optimization of selected proofs. The authors are responsible for all mathematical results, references, and the final version of the manuscript.

\end{document}